\documentclass[12pt]{amsart}
\usepackage{amsmath}
\usepackage{amssymb}

\def\underTilde#1{{\baselineskip=0pt\vtop{\hbox{$#1$}\hbox{$\sim$}}}{}}
\def\undertilde#1{{\baselineskip=0pt\vtop
  {\hbox{$#1$}\hbox{$\scriptscriptstyle\sim$}}}{}}

 \newtheorem{thm}{Theorem}
 
 \newtheorem{lem}[thm]{Lemma}

 \theoremstyle{definition}

 \theoremstyle{remark}

\newtheorem{claim}[thm]{Claim}

\newcommand{\A}{ A_{\mathsf{W}} }
 \numberwithin{equation}{section}

\begin{document}

\title{Bounded Martin's Maximum With Many Witnesses}

\author{Stuart Zoble \\ December 6, 2008}

\address{Department of Mathematics, Wesleyan University}
\email{azoble@wesleyan.edu}



\begin{abstract} We study a strengthening of Bounded Martin's Maximum which asserts
that if a $\Sigma_1$ fact holds of $\omega_2^V$ in a stationary set preserving extension then it holds in $V$ for a stationary set of ordinals less than $\omega_2$.  We show that this principle implies Global Projective Determinacy, and therefore does not hold in the $\mathbb{P}_{max}$ model for $\mathsf{BMM}$, but that the restriction of this principle to forcings which render $\omega_2^V$ countably cofinal does hold in the $\mathsf{BMM}$ model, though it is not a consequence of $\mathsf{BMM}$.    \end{abstract}

\maketitle

\section*{Introduction}  \noindent Bounded Martin's Maximum, denoted $\mathsf{BMM}$, is the assertion that $$(H(\omega_2),\in) \prec_{\Sigma_1} (H(\omega_2),\in)^{V^{\mathbb{P}}}$$ \noindent whenever $\mathbb{P}$ is
stationary set preserving.  Because $\Sigma_1$ facts are upward absolute and $Col(\omega_1,\omega_2)$ can be appended to a given stationary set preserving forcing, the formulation below is equivalent.

\begin{quote} $\mathsf{BMM}$ denotes the following assertion.  Suppose $a \in
H(\omega_2)$, $\varphi(x,a)$ is a $\Sigma_1$ formula, $\mathbb{P}$ is
stationary set preserving forcing notion and whenever $G \subset \mathbb{P}$ is
generic
$$(H(\omega_2),\in)^{V[G]} \models \varphi(\omega_2^V,a).$$ \noindent
Then $$(H(\omega_2),\in)^{V} \models \varphi(\delta,a)$$ \noindent for some
ordinal $\delta < \omega_2$.  \end{quote}

\noindent In this paper we will study a strengthening of $\mathsf{BMM}$ which asserts that if a $\Sigma_1$ fact holds of $\omega_2^V$ in a stationary set preserving
extension then it holds in $V$ for a stationary set of ordinals less than $\omega_2$.  In other words, we will
strengthen $\mathsf{BMM}$ as above by replacing the phrase "some ordinal $\delta < \omega_2$" with "a stationary set of ordinals $\delta < \omega_2$".

\begin{quote} $\mathsf{BMM^{s}}$ denotes the following assertion.  Suppose $a \in
H(\omega_2)$, $\varphi(x,a)$ is a $\Sigma_1$ formula, $\mathbb{P}$ is
stationary set preserving forcing notion and whenever $G \subset \mathbb{P}$ is
generic
$$(H(\omega_2),\in)^{V[G]} \models \varphi(\omega_2^V,a).$$ \noindent
Then $$(H(\omega_2),\in)^{V} \models \varphi(\delta,a)$$ \noindent for a
stationary set of ordinals $\delta < \omega_2$.  \end{quote}

\noindent A definition of $\mathsf{BMM^{s,++}}$ is obtained by replacing the occurrences of the structure $(H(\omega_2),\in)$ with the structure $(H(\omega_2),\in,\mbox{NS})$ in the definition above, and regarding $\varphi$ as $\Sigma_1$ in the expanded language.

Our main theorem is that $\mathsf{BMM^{s,++}}$ implies
that Projective Determinacy holds in all set-generic extensions.  To put this in perspective, it is open whether $\mathsf{BMM}$ implies $\underTilde{\Delta}^1_2$ determinacy, and known that $\mathsf{BMM}$ does not imply $\underTilde{\Delta}^1_3$ determinacy in the universe after $\omega_2$ is collapsed as Theorem 10.99 of \cite{W} shows that $$L^{M_{1}^{\#}}(\mathbb{R})[G] \models \mathsf{BMM}^{++}$$ \noindent where $N = L^{M_{1}^{\#}}(\mathbb{R})$ is the minimal inner model closed under the $M_{1}^{\#}$ operation and containing $\mathbb{R}$, $N$ satisfies $\mathsf{AD}$, and
$G \subset \mathbb{P}_{max}$ is $N$-generic.  Thus, it seems that $\mathsf{BMM^{s,++}}$ is a bit stronger than $\mathsf{BMM}$ and in particular does not hold in the model described above.  We are not able to prove $\mathsf{PD}$ from $\mathsf{BMM^s}$ alone, but we can if we assume in addition that the nonstationary ideal is saturated, and the proof produces for example a set of ordinals $E$ for which $M_{2}^{\#}(E)$ exists and
$$P(\omega_1) \subset M_2^{\#}(E),$$ \noindent which is enough to conclude that $\mathsf{BMM^s}$ fails in the $\mathsf{BMM}$ model.  We are able to show, however, that a special case of this forcing axiom does hold in the $\mathsf{BMM}$ model, namely
$\mathsf{BMM^{s_0}}$, which will denote the restriction of $\mathsf{BMM^s}$ to stationary set preserving forcings $\mathbb{P}$ for which $$V[G] \models cf(\omega_2^V) = \omega$$ \noindent whenever $G \subset \mathbb{P}$ is $V$ generic.

\begin{thm} \

\begin{enumerate}

\item $\mathsf{BMM^{s,++}}$ implies $\mathsf{PD}$ in all generic extensions.
\item $\mathsf{BMM^s}$ fails in the $\mathsf{BMM}$ model.

\item $\mathsf{BMM^{s_0,++}}$ holds in the $\mathsf{BMM}$ model.

\end{enumerate}
\end{thm}

 The proof of the theorem makes essential use of the well ordering of $P(\omega_1)$ given
under $\mathsf{BPFA}$ by Caicedo and Velickovic in \cite{CV}.  Once one understands why their
well-ordering is $\underTilde{\Delta}_1$ over $H(\omega_2)$ it is easy to construct from it a bijection $$\mathsf{W}:\omega_2 \rightarrow P(\omega_1)$$ \noindent whose initial segments
are uniformly $\underTilde{\Sigma}_1$ definable over $H(\omega_2)$.  That is, $\mathsf{W}$ is a bijection and
there is a $\Sigma_1$ formula $\psi$ and a parameter $a \in H(\omega_2)$ such that for every $x \in H(\omega_2)$ and $\beta < \omega_2$,
$$(H(\omega_2),\in ) \models \psi(x,\beta,a) \ \mbox{
iff } x = \mathsf{W} \upharpoonright \beta.$$

\noindent Well-orderings due to Moore and Todorcevic do not seem to suit our purposes, though Woodin's well-ordering from $\psi_{AC}$ does give such a set $\mathsf{W}$ which is uniformly $\underTilde{\Sigma}_1$ definable over the structure $(H(\omega_2),\in,NS)$.  The reader curious about how the well-ordering is used to increase the expressive
power of the $\Sigma_1$ formula in the definition of $\mathsf{BMM^s}$ could skip directly
to Lemma 9 below.

We now give some further background information.  $\mathsf{BMM}$ implies
a bounded version of the strong reflection principle which we denote by $\mathsf{BSRP}(\omega_2)$,
and which asserts that any projective stationary subset $S$ of
$[\omega_2]^{\omega}$ which is $\underTilde{\Sigma}_1$-definable over the structure $(H(\omega_2),\in)$ reflects to a club in $[\gamma]^{\omega}$ for some, equivalently unboundedly many, ordinals $\gamma < \omega_2$.  The nucleus of this paper, now essentially the base case of the $\mathsf{PD}$ induction, was the observation that some open questions regarding $\mathsf{BMM}$ could be solved assuming in addition the following enhanced version of $\mathsf{BSRP}(\omega_2)$ giving a stationary set of club reflection points.  The axiom $\mathsf{BMM^s}$ is the natural generalization of $\mathsf{BMM}$ which implies this stronger reflection principle.

\begin{quote} $\mathsf{BSRP^{s}}(\omega_2)$ denotes the assertion that any projective stationary set
$S \subset [\omega_2]^{\omega}$ which is $\Sigma_1$-definable in
$H(\omega_2)$ reflects to a club stationarily often in the sense that $S \cap [\gamma]^{\omega}$ contains a club in
$[\gamma]^{\omega}$ for a stationary set of $\gamma < \omega_2$.  \end{quote}

The particular questions which interest us ask whether certain consequences of Martin's Maximum in fact follow from $\mathsf{BMM}$, for example those on the following list.

\begin{enumerate}

\item The nonstationary ideal is precipitous

\item Woodin's principle $\psi_{AC}$

\item Canonical function bounding

\item $\undertilde{\delta}^1_2 = \omega_2$

\item $\underTilde{\Delta}^1_2$ determinacy

\end{enumerate}

Of course, the determinacy question is really just a question of consistency strength.  The best result
to date is due to Schindler who shows in \cite{Sc2} that $\mathsf{BMM}$ implies the existence of
an inner model with a strong cardinal.  Regarding $\undertilde{\delta}^1_2 = \omega_2$, Woodin
gets this from NS saturated with a measurable cardinal in \cite{W}, and hence from
$\mathsf{BMM}$ together with a Woodin cardinal and a measurable above using a theorem of Shelah.  Schindler
and Claverie have recently proved $\undertilde{\delta}^1_2 = \omega_2$ from $\mathsf{BMM}$ together
together with NS precipitous.  Canonical function bounding follows outright from $\psi_{AC}$ by an argument of Aspero in \cite{AW}, and Woodin obtains $\psi_{AC}$ as a consequence of $\mathsf{BMM}$ with either
a measurable cardinal or NS precipitous assumed in addition (see \cite{W}).  Our initial observation
went as follows.

\begin{thm} Assume $\mathsf{BMM}$ and $\mathsf{BSRP^{s}}(\omega_2)$.  Then $\underTilde{\Delta}^1_2$ determinacy, $\undertilde{\delta}^1_2 = \omega_2$, and canonical function bounding holds.
\end{thm}

\noindent The theorem is proved by showing that $$\mathsf{W}^{\dagger} \models \mbox{NS saturated}$$ \noindent where NS denotes the nonstationary ideal on $\omega_1$.  This suggests a new entry for the above list of possible consequences of $\mathsf{BMM}$.

\begin{quote} Does $\mathsf{BMM^{++}}$ imply NS saturated in $L(P(\omega_1))$?
\end{quote}

As one would suspect, $\mathsf{BSRP^s}(\omega_2)$ is not a consequence of $\mathsf{BMM}$, and
we establish this by way of the Tilde operation.  Recall that $\tilde{T}$,
for a subset $T$ of $\omega_1$, is defined to be the set of $\alpha < \omega_2$ for which there is a club of
$\sigma \in [\alpha]^{\omega}$ whose order type belongs to $T$.  We show that in a forcing
extension of a model satisfying Martin's Maximum,
$\mathsf{BMM}$ holds and the nonstationary ideal is saturated, yet there exists a stationary set $T \subset \omega_1$ for which $\tilde{T}$ is nonstationary in $\omega_2$.  An argument of Larson
from \cite{L} shows that such a set $\tilde{T}$ must be stationary under $\mathsf{BSRP^s}(\omega_2)$ together
with NS saturated, and so $\mathsf{BSRP^s}(\omega_2)$ must fail in this model.

\begin{thm} $\mathsf{BSRP^s}(\omega_2)$ is not a consequence of $\mathsf{BMM^{++}}$ together with the saturation
of the nonstationary ideal.
\end{thm}

In a similar vein, arguments of Larson from \cite{L} coupled with a Theorem of Woodin from \cite{W} will produce
models in which $\mathsf{BMM^{s_0}}$ fails but $\mathsf{BMM}$ as well as other hypotheses hold.  Finally, we give an application of $\mathsf{BMM^s}$ which does not seem to have anything to
do with stationary reflection but involves rather the notion of a {\em disjoint club sequence on} $\omega_2$, an invention of Krueger from \cite{K} who derives one from $\mathsf{MM}(c)$.  We observe here that $\mathsf{BMM^s}$ implies the existence of such a sequence, and this will allow us to separate $\mathsf{BMM}$ and some
consequences of $\mathsf{BMM^s}$ used in the proof of Theorem 1, from $\mathsf{BMM^s}$ itself.

This paper is organized as follows.  We start with a brief discussion of the $\Sigma_1$ well-ordering and
 then prove the two results (Theorems 2 and 3) concerning $\mathsf{BSRP^s}(\omega_2)$.  We then give the
$\mathsf{PD}$ proofs, followed by the $\mathbb{P}_{max}$ argument, and close with the other separation result.
We would like to thank Paul Larson for directing us to the relevant results in \cite{L}, and for many
enlightening conversations.

\section{Results} We need that the wellordering of \cite{CV} gives rise to a
uniformly $\underTilde{\Sigma}_1$ enumeration of $P(\omega_1)$ as described in the introduction.  For
the reader's convenience we describe how this is obtained.

\begin{lem} (Caicedo, Velickovic \cite{CV})  Assume $\mathsf{BPFA}$.  Then there is a bijection $$\mathsf{W}:\omega_2 \rightarrow P(\omega_1)$$ \noindent whose initial segments are uniformly $\underTilde{\Sigma}_1$ definable over
$H(\omega_2)$.  \end{lem}

\begin{proof}  The parameter involved in the
definition is a certain subset $c \subset \omega_1$.  This parameter gives rise to an $\omega_1$-sequence $d$
of pairwise almost disjoint elements of $[\omega]^{\omega}$ which will be used as well.  The authors of \cite{CV} define
a notion of a triple $\alpha < \beta < \gamma < \omega_2$ coding a real $r$ which is a
$\Sigma_1$ notion in the parameter $c$.  Fixing a reasonable way of coding a triple by a single ordinal
and composing, we thus have a $\Sigma_1$ formula which says that an ordinal $\delta < \omega_2$ codes a real $r$, and they prove that every real is so coded.  Every subset $a$ of $\omega_1$
is coded by a real $r \subset \omega_1$ via $d$ in the standard way since $\mathsf{MA}_{\omega_1}$ holds.  Let $T$ be the theory described in \cite{CV} which includes the sentence asserting that every real is coded by an ordinal, as well as $\mathsf{MA}_{\omega_1}$, among other axioms.  Then $H(\omega_2) \models T$ and for transitive models $M,N$ of $T$

\begin{enumerate}

\item If $M \cap \omega_2 = N \cap \omega_2$ then $M = N$

\item If $M \cap \omega_2 < N \cap \omega_2$ then $M \in N$.

\end{enumerate}

\noindent  For a real $r$ let $N_{r}$ be the least model of $T$ with $r \in N_{r}$.  Then a well-ordering of the reals is given by $r \prec s$ if $N_{r} \in N_{s}$ or $N_{r} = N_{s} = N$ and the least ordinal which
codes $r$ in the sense of $N$ is less than the least ordinal which
codes $s$ in the sense of $N$.  Let $\{ r_{\delta} \ | \ \delta < \omega_2 \}$ be the enumeration
of the reals according to this ordering, and set $\mathsf{W}_{0}(\delta) = a$ if $r_{\delta}$ codes
$a$ via $d$.  Then $\mathsf{W}_{0}$ is a uniformly $\underTilde{\Sigma_1}$ definable surjection therefore
gives rise to such a uniformly definable bijection $\mathsf{W}$ in the obvious way.  \end{proof}

\noindent We now show how $\mathsf{BSRP^{s}}(\omega_2)$ can be used to
prove that NS is saturated in an inner model with a measurable cardinal which
contains $P(\omega_1)$.  We use Schindler's theorem from \cite{Sc2} to produce the
measurable, although this can be avoided if $\mathsf{BPFA}^{++}$ is assumed in place
of $\mathsf{BMM}$.

\begin{thm}  Assume either $\mathsf{BMM}$ or $\mathsf{BPFA}^{++}$.  Assume that
$\mathsf{BSRP^{s}}(\omega_2)$ holds in addition.  Then

\begin{enumerate}

\item $\underTilde{\Delta}^1_2$ determinacy holds

\item Every function from $\omega_1$ to $\omega_1$ is bounded by a canonical function

\item $\undertilde{\delta}^1_2 = \omega_2$.

\end{enumerate}

\end{thm}

\begin{proof} Since $\mathsf{BPFA}$ holds we can let $\mathsf{W}$ denote the
unifomly $\Sigma_1$ enumeration of $P(\omega_1)$ given by Lemma 4.  For convenience
we will think of $\mathsf{W}$ as a subset of $\omega_2 \times \omega_1$ with the property that
$$P(\omega_1) = \{ \mathsf{W}_\alpha \ | \ \alpha < \omega_2 \}$$ \noindent where
$\mathsf{W}_{\alpha}$ denotes the set $\{ \gamma \ | \ (\alpha,\gamma) \in \mathsf{W} \}$.  Thus
there is a $\Sigma_1$ formula $\phi(x,y,z)$ and a parameter $a \in H(\omega_2)$ such that
for every $x \in H(\omega_2)$ and $\beta < \omega_2$
$$(H(\omega_2),\in ) \models \psi(x,\beta,a) \ \mbox{
iff } x = \mathsf{W} \cap (\beta \times \omega_1).$$

\noindent First assume $\mathsf{BMM}$ and $\mathsf{BSRP^{s}}(\omega_2)$.
Schindler has shown (see \cite{Sc1} and \cite{Sc2}) that $X^{\dagger}$ exists for every set $X$ under $\mathsf{BMM}$ so in particular $\mathsf{W}^{\dagger}$ exists.  For any set $X$ we let $\mathcal{M}(X)$ denote $X^{\dagger}$.
We will use $\mathsf{BSRP^{s}}(\omega_2)$ to seal a putative least bad antichain in $\mathcal{M}(\mathsf{W})$ thereby showing that
$$\mathcal{M}(\mathsf{W}) \models \mbox{ NS is saturated}.$$
\noindent Let us assume toward a contradiction that in $\mathcal{M}(\mathsf{W})$
there is a maximal antichain in $P(\omega_1)/NS$ of
size $\omega_{2}^{\mathcal{M}(\mathsf{W})} = \omega_2$.  Let $\mathcal{A}$ denote the least antichain
in the canonical well-ordering of $\mathcal{M}(\mathsf{W})$.  Using $\mathsf{W}$ we may
code $\mathcal{A}$ as a subset $A$ of $\omega_2$ given by
$$A = \{\alpha < \omega_2 \ | \ \mathsf{W}_{\alpha} \in \mathcal{A} \}.$$
\noindent For $\sigma \in [\omega_2]^{\omega}$ let
$\pi_{\sigma}:\sigma \rightarrow otp(\sigma)$ be the collapse of
$\sigma$, let $\mathsf{W}_{\sigma}$ denote the image $\pi_{\sigma}[\mathsf{W} \cap
\sigma]$, noting that $\pi_{\sigma}$ acts on pairs in the obvious way, and let $A_{\sigma} = \pi_{\sigma}[A \cap \sigma]$.  For a
club of $\sigma$ it will be true that the code of the least
antichain of $\mathcal{M}(\mathsf{W}_{\sigma})$ is $A_{\sigma}$.  In every
case, let us use $A_{\sigma}$ to denote the code of the least maximal antichain of length $\omega_2$ in the sense of
$\mathcal{M}(\mathsf{W}_{\sigma})$ if it exists.  Define the set $S \subset
[\omega_{2}]^{\omega}$ to consist of all $\sigma$ satisfying

\begin{enumerate}

\item $\mathcal{M}(\mathsf{W}_{\sigma})$ thinks that $\mathsf{W}_{\sigma} \subset \omega_2 \times \omega_1$
enumerates $P(\omega_1)$ in length $\omega_2$

\item $\mathcal{M}(\mathsf{W}_{\sigma})$ thinks that the least $NS$
antichain exists and is coded as above by some $A_{\sigma} \subset
otp(\sigma)$

\item There exists $\alpha \in \sigma$ such that
$\pi_{\sigma}(\alpha) \in A_{\sigma}$ and $\sigma \cap \omega_1 \in
\mathsf{W}_{\alpha}.$

\end{enumerate}

\noindent This set is $\underTilde{\Sigma}_1$ definable over
$(H(\omega_2),\in)$.  To verify $\sigma \in S$ it suffices to find
an ordinal $\delta > \omega_1$ with $sup(\sigma) \subset \delta$ and a transitive set $N$
with $\delta \subset N$ which satisfies enough set theory, computes
$\mathsf{W} \cap \delta \times \omega$ correctly, contains $\sigma$, contains a
$\mathsf{W}_{\sigma}$ premouse $M$, thinks that $M = \mathcal{M}(\mathsf{W}_{\sigma})$
and that the conditions above are satisfied.  Note that any such structure is correct about
$M = \mathcal{M}(\mathsf{W}_{\sigma})$ since it is a $\Pi^1_2$ condition.  We claim that $S$ is projective stationary.  Since $P(\omega_1) \subset \mathcal{M}(\mathsf{W})$ the
antichain coded by $A$ is truly a maximal antichain in
$P(\omega_1)/NS$.  It is well known that the set of $\sigma \in
[\omega_2]^{\omega}$ such that $$\sigma \cap \omega_1 \in
\bigcup_{\alpha \in \sigma \cap A} \mathsf{W}_{\alpha}$$ \noindent is
projective stationary.  Our set $S$ differs from this set on a
nonstationary subset of $[\omega_2]^{\omega}$.  Let $\theta$ be large enough and
let $(H_{\xi} \ | \ \xi < \omega_2)$ be an increasing sequence of elementary submodels
of $H(\theta)$, each of size $\omega_1$, so that $\mathcal{M}(\mathsf{W}) \in H_{\xi}$ and $H_{\xi} \cap \omega_2 \in \omega_2$ for each $\xi$.  We may assume that the sequence $H_{\xi} \cap \omega_2$ is strictly increasing, continuous, and converges to $\omega_2$.  Let $C \subset \omega_2$ be a club so that $H_{\xi} \cap \omega_2 = \xi$ for
$\xi \in C$.  Thus there exists a $\delta <
\omega_2$ such that

\begin{enumerate}

\item[(a)] $\delta \in C$ and

\item[(b)] $S$ contains a club in $[\delta]^{\omega}$.

\end{enumerate}

\noindent  Let $\pi:H_{\delta} \rightarrow H$ be the transitive collapse.  Thus
$$\pi( \mathcal{M}(\mathsf{W})) = \mathcal{M}(\mathsf{W} \cap (\delta \times \omega_1))$$ \noindent so that
the least antichain of $\mathcal{M}(\mathsf{W} \cap (\delta \times \omega_1))$ is coded by $\pi(A) = A \cap \delta$.
Now we may let $(N_{\gamma} \ | \ \gamma < \omega_1)$ be an increasing sequence of countable, elementary submodels of
$H(\theta)$ which contain $\delta$ so that $(\sigma_{\gamma} \ | \ \gamma < \omega_1)$ is continuous
and exhaustive in $[\delta]^{\omega}$, where $\sigma_{\gamma} = N_{\gamma} \cap \delta$.  It follows by (b) that
there is a club $D \subset \omega_1$ so that $\sigma_{\gamma} \in S$ for $\gamma \in D$.  This implies that
the diagonal union of the sets coded in $A \cap \delta$ contains $D$, which is a
contradiction.  Thus in $\mathcal{M}(\mathsf{W})$ there is a measurable
cardinal and $NS$ is saturated.  It follows $\underTilde{\Delta}^1_2$-determinacy and $\undertilde{\delta}^1_2 = \omega_2$
hold in $\mathcal{M}(\mathsf{W})$ by 7.1 of \cite{S} and 3.17 of \cite{W} respectively.
They therefore hold in $V$ as $P(\omega_1) \subset \mathcal{M}(\mathsf{W})$.  Canonical function bounding is a consequence of NS saturated so it holds in $\mathcal{M}(\mathsf{W})$ and hence in $V$.   Now we assume
$\mathsf{BPFA}^{++}$ in place of $\mathsf{BMM}$.  Thus, we have to do without Schindler's theorem.  The argument
above shows that $$L[\mathsf{W}] \models \mbox{NS is saturated}$$ \noindent using only $\mathsf{BPFA}$ and
$\mathsf{BSRP}^{s}(\omega_2)$.  This is because if $L[\mathsf{W}]$ thinks that there is a maximal antichain
$\mathcal{A}$ of size $\omega_2$ then we can define an operation $\mathcal{M}(x)$ which associates
to a set of ordinals $x$ the least level of $L[x]$ which satisfies a sufficient fragment of
set theory, thinks that $x \subset \omega_2 \times \omega_1$ codes $P(\omega_1)$ and that
such an antichain exists.  Virtually the same argument yields a contradiction.  Working inside
$L[W]$ we can argue that $x^{\dagger}$ exists for every bounded subset $x$ of $\omega_1$ and
use the generic ultrapower to conclude that $P(\omega_1)$ is closed under daggers.  We may now
simply quote 10.108 of \cite{W} to conclude that $\mathsf{W}^{\dagger}$ exists, and in fact that every
set has a dagger, but we will elaborate on this point as a generalization of this argument will play a role
in the $\mathsf{PD}$ proof.  Assume toward a contradiction that there is a set of ordinals $A$ so that $A^{\dagger}$ does not exist.  Let $\theta$
be a cardinal containing $A$ and let $g \subset Col(\omega_1,\theta)$ be $V$-generic.  In $V[g]$ let
$a \subset \omega_1$ be such that $$A \in L[a].$$ \noindent By standard arguments we must
have that $$V[g] \models a^{\dagger} \mbox{ does not exist},$$
\noindent but that $$V[g] \models (a \cap \alpha)^{\dagger} \mbox{ does exist}$$ \noindent for every $\alpha < \omega_1$ as no new reals were added.  Thus in $V[g]$ there are stationary
sets $S,T$ and $p$ such that $\alpha \in S$ implies $p \in (a \cap \alpha)^{\dagger}$ and $\alpha \in T$ implies
$p \notin (a \cap \alpha)^{\dagger}$.  This can be reformulated as a $\underTilde{\Sigma}_1$ fact in the language
of set theory together with a predicate for NS.  Since $Col(\omega_1,\theta)$ is proper there are such sets $\bar{a},\bar{S},\bar{T},\bar{p}$ in $V$ which
implies that $\bar{a}^{\dagger}$ cannot exist in $V$, a contradiction.  Another application of the sealing argument gives the saturation of NS in $\mathsf{W}^{\dagger}$ from which the result follows.  \end{proof}

From $\mathsf{BSRP^{s}}(\omega_2)$ together with a uniformly $\Sigma_1$ enumeration of
$P(\omega_1)$ the argument above shows that $L(P(\omega_1))$ is a model of $\mathsf{ZFC}$ together with
NS saturated, and hence bounding holds and there is an inner model with a Woodin cardinal by a recent result from \cite{SJ}.  We conjecture that $\mathsf{BSRP^{s}}(\omega_2)$ alone implies an inner model with a Woodin cardinal.
Before moving on to global principles we demonstrate that $\mathsf{BSRP^{s}}(\omega_2)$ is not a
consequence of $\mathsf{BMM}$.  Recall that $\tilde{T}$ is the set of $\alpha < \omega_2$ for which there is
a club of $\sigma \in [\alpha]^{\omega}$ with $otp(\alpha) \in T$.  Paul Larson pointed out to us
that in the presence of NS saturated, the principle
$\mathsf{BSRP^s}(\omega_2)$ would imply that $\tilde{S}$ is stationary for every stationary
set $S \subset \omega_1$ (see 3.9 of \cite{L}), and that a model of his from \cite{L} could be modified
to produce a model of $\mathsf{BMM}^{++}$ together with a stationary set $S$ so that
$$\{ \alpha \in \tilde{S} \ | \  cf(\alpha) = \omega \}$$ \noindent is nonstationary.  A subtle point, however, is that $\mathsf{MM}(c)$ holds in this model, and so $\mathsf{BSRP^s}(\omega_2)$ doesn't always give a stationary set of cofinality $\omega$ reflection points.  It turns out that, assuming $\mathsf{MM}$, we can shoot a club through some
$\tilde{T}$ without adding subsets of $\omega_1$ while preserving the
saturation of NS.  This will be enough to separate $\mathsf{BMM}$ from $\mathsf{BSRP^s}(\omega_2)$.

\begin{lem} (Larson) Assume $\mathsf{BSRP^s}(\omega_2)$ holds and that the nonstationary ideal
is saturated.  Then $\tilde{T} \subset \omega_2$ is stationary for every stationary
set $T \subset \omega_1$.  \end{lem}

\begin{proof}  Fix stationary sets $S,T \subset \omega_1$.  Fix a function $f:[\omega_2]^{\omega} \rightarrow \omega_2$ and let $j: V \rightarrow M \subset V[G]$ be the NS generic ultrapower.  Let $H \subset j(P(\omega_1)/NS)$ be $M$-generic with ultrapower map $k:M \rightarrow N$.  We
assume that $S \in G$ and $j(T) \in H$.  Let $\sigma$ denote $(k \circ j)[\omega_2^V]$.  Then $\sigma$ is
a countable subset of $(k \circ j)(\omega_2)$ with the property that $\sigma \cap \omega_1^N \in k \circ j(S)$ and
$otp(\sigma) \in  k \circ j(S)$.  Further $\sigma$ is closed under $(k \circ j)(f)$.  $N$ must see such a countable set with these properties and so by
elementarity and the fact that $S$ was arbitrary we conclude that $$\{ \sigma \in [\omega_2]^{\omega} \ | \ otp(\sigma) \in T \}$$ \noindent
is projective stationary in $V$.  Since this set is $\Sigma_1$ definable from $T$ as a parameter,
we get the desired conclusion.  \end{proof}

\begin{thm}  Assume $\mathsf{MM}$.  Then there is a forcing notion $\mathbb{P}$ of size $\omega_2$
such that whenever $G \subset \mathbb{P}$ is $V$-generic then $V[G]$ satisfies

\begin{enumerate}

\item $\mathsf{BMM}^{++}$,

\item the nonstationary ideal is saturated, and

\item there is a stationary set $T \subset \omega_1$ with $\tilde{T}$ nonstationary.

\end{enumerate}

\noindent Thus $\mathsf{BSRP^s}(\omega_2)$ must fail in $V[G]$.
 \end{thm}

\begin{proof}  Fix a stationary and costationary set $T \subset \omega_1$ and assume $\mathsf{MM}$.  Let $\mathbb{P}$
be the poset for shooting a club through $\tilde{T}$.  $\mathbb{P}$ consists of closed subsets
of $\tilde{T}$ of size $\omega_1$ ordered by end extension.  We first show that $\mathbb{P}$ does not add new subsets of $\omega_1$.  Let $\mathcal{S}$ denote $$\{ \sigma \in [\omega_2]^{\omega} \ | \ \mbox{sup}(\sigma) \in \tilde{T} \mbox{ and } \mbox{otp}(\sigma) \in T \}.$$

\begin{claim} $\mathcal{S}$ is projective stationary.
\end{claim}

\begin{proof}  Fix a stationary set $S \subset \omega_1$. Let $G \subset P(\omega_1)/NS$ be $V$-generic
with $S \in G$.  Let $\pi:V \rightarrow M \subset V[G]$ be the generic embedding with $G \subset P(\omega_1)/NS$.
Let $f:[\omega_2]^{<\omega} \rightarrow \omega$ be arbitrary.  Let $C$ denote the set of $\delta < \omega_2$
with $f[[\delta]^{<\omega}] \subset \delta$.  By Lemma 5.8 of \cite{W},
$$\{ \alpha \in \tilde{T} \ | \  cf(\alpha) = \omega \}$$ \noindent is a stationary subset of $\omega_2$,
and so there is $\delta \in \tilde{T} \cap C$ such that $\delta$ has countable cofinality.  Let $\sigma = \pi[\delta]$.  Then $\sigma$ has the following properties in $V[G]$.

\begin{enumerate}

\item $\delta = \mbox{otp}(\sigma) \in \pi(T)$

\item $\mbox{sup}(\sigma) = \pi(\delta)$ belongs to the tilde of $\pi(T)$

\item $\sigma \cap \omega_1 \in \pi(S)$.

\item $\sigma$ is closed under $\pi(f)$

\end{enumerate}

\noindent Since $\sigma$ is countable in $V[G]$ we have $\sigma \in M$ and so by elementarity
we find in $V$ such a set $\sigma$ in $\mathcal{S}$ which is closed under $f$ and with $\sigma \cap \omega_1 \in S$ as desired.  \end{proof}

\noindent Let $\theta$ be large and let $\mathcal{S}^{*}$ denote the set of countable elementary
submodels $X \prec H(\theta)$ with $\mathbb{P} \in X$ and $X \cap \omega_2 \in \mathcal{S}$.  Let $\tau$
be a term for a subset of $\omega_1$.  Using $\mathsf{MM}$ let $$(X_{\alpha} \ | \ \alpha < \omega_1)$$ \noindent be
an increasing, continuous $\in$ chain with each $$X_{\alpha} \cap \omega_2 \in \mathcal{S}^{*},$$ \noindent and $\tau \in X_{0}$.  Inductively define a decreasing sequence of conditions $p_{\alpha}$ so that $p_{\alpha} \in X_{\alpha+1}$, $$p_{\alpha} = q \cup \{ \mbox{sup}(X_{\alpha} \cap \omega_2) \}$$ \noindent where $q$ is an $X$-generic
filter.  We can assume that $$\bigcup_{\alpha < \omega_1} X_{\alpha} \cap \omega_2 = \delta$$ \noindent for
some $\delta < \omega_2$ and so $\delta \in \tilde{T}$.  It follows that
$$p = (\bigcup_{\alpha < \omega_1} p_{\alpha}) \cup \{ \delta \}$$ \noindent decides $\tau$.  Thus $\mathbb{P}$
does not add new subsets of $\omega_1$.  Now let $\tau$ be a term and $p_0$ a condition such that $p$
forces that $\tau$ is a function from $\omega_2$ to $P(\omega_1)/NS$ whose range is
a maximal antichain.  We suppress $p_0$.  Let $\mathcal{S}^{**}$ consist of all $X \in [H(\omega_2)]^{\omega}$ satisfying

\begin{enumerate}

\item $X \prec H(\omega_2)$

\item $X \cap \omega_2 \in \mathcal{S}^{*}$

\item For a dense set of $q \in \mathbb{P} \cap X$ there is an ordinal $\gamma \in X$ so
that $$q \Vdash_{\mathbb{P}}^{V} X \cap \omega_1 \in \tau(\gamma).$$

\end{enumerate}

\noindent We claim that $\mathcal{S}^{**}$ is projective stationary.  Otherwise by pressing
down we find a condition $q_0$ and a stationary set $\mathcal{T}$ of $X \prec H(\omega_2)$
with $X \cap \omega_2 \in \mathcal{S}^{*}$ such that for all $q \in \mathbb{P} \cap X$ with
$q \leq q_0$ and all $\gamma \in X$, $$\neg \ (q \Vdash_{\mathbb{P}}^{V} X \cap \omega_1 \in \tau(\gamma)).$$
\noindent Since NS is saturated the set $\mathcal{T}$ must be $A$-projective stationary for
some stationary set $A \subset \omega_1$.  Let $q_1$ be a condition below $q$, $\gamma < \omega_2$
and $B \subset A$ so that $q_1$ forces $\tau(\gamma) \cap A = B$.  There are stationary many
$X \in \mathcal{T}$ such that $X \cap \omega_1 \in B$, $\gamma \in X$, and $q_1 \in X$.  Any such
$X$ gives the desired contradiction.  Now let $G \subset \mathbb{P}$ be $V$-generic.  We have
shown that $\mathsf{BMM}^{++}$ holds and NS is saturated in $V[G]$.  Moreover, $\tilde{T}$ contains a club so that $$\tilde{S} \subset \omega_2 \setminus \tilde{T}$$ \noindent must be nonstationary where $S =
\omega_1 \setminus T$.  Since $S$ is stationary, $\mathsf{BSRP^s}(\omega_2)$ must fail by Lemma 6.  \end{proof}

The $\mathsf{BPFA}^{++}$ part of the $\underTilde{\Delta}^1_2$ determinacy proof illustrates the two main elements, {\em sealing} and {\em lifting},
of the proof of $\mathsf{PD}$ from $\mathsf{BMM^{s,++}}$.  This proof is modeled on Woodin's proof of $\mathsf{PD}$ from $\mathsf{MM}(c)$ in which the saturation of NS is the key hypothesis used in contradicting the existence of the core model.  While it is unlikely that $\mathsf{BMM^{s,++}}$ implies the saturation of NS, it will imply saturation inside the various inner models of interest as we proceed through the $\mathsf{PD}$ induction.  The other element of the $\mathsf{MM}(c)$ proof involves lifting closure under fine structural operations from $P(\omega_1)$ to $P(\omega_2)$ using simultaneous reflection as in 9.78 of \cite{W}.  With $\mathsf{BMM^{s,++}}$ we can lift closure from $P(\omega_1)$ to $V$ by an argument which is more in the spirit of 10.108 of \cite{W}, and this is where the "++" seems unavoidable.  We can however get by with $\mathsf{BMM^s}$ in this connection if NS is saturated in addition, though we can only lift closure from $P(\omega_1)$ to definable subsets of $\omega_2$.  This is enough to get $\mathsf{PD}$ and is how we will show that $\mathsf{BMM^s}$ fails in the $\mathsf{BMM}$ model.

We first derive the following definable version of $\mathsf{MM}(c)$ to illustrate how
the definable wellordering is used to increase the expressive power of the $\Sigma_1$ formula 
appearing in the definition of $\mathsf{BMM^s}$.  As a corollary we get a version of $\mathsf{BSRP^s}(\omega_2)$
for all definable projective stationary sets which will be used in Theorem 11 below.

\begin{lem} Assume $\mathsf{BMM^s}$.  Suppose  Suppose $\mathbb{P}$ and
$D$ are first order definable over $H(\omega_2,\in)$, $\mathbb{P}$ is a poset, and $D$ is a partial map
from $H(\omega_2)$ to $H(\omega_2)$ with the property that $D(a) \subseteq \mathbb{P}$ is dense where
defined.  Then there is a stationary set of $\delta < H(\omega_2)$ such that there exists $X$ and $G$
satisfying

\begin{enumerate}

\item $X$ is a transitive and fully elementary submodel of $H(\omega_2)$

\item $X \cap \omega_2 = \delta$

\item $G$ is a filter on $\mathbb{P} \cap X$

\item $D(a) \cap G \neq \emptyset$ \noindent for any $a \in X \cap \mbox{dom}(D)$.

\end{enumerate}

\end{lem}

\begin{proof}  Let $\mathbb{P}$ and $D$ be as above.  Let
$G \subset \mathbb{P}$ be $V$-generic.  Let $\psi(x,y)$ be the formula
(with parameter suppressed) which defines the initial segments of $\mathsf{W}$ uniformly over $H(\omega_2)$.
Note that for any $\beta < \omega_2^V$ the set $\mathsf{W} \upharpoonright \beta$ is the unique
witness to $\psi(x,\beta)$ in $V$ as well as $V[G]$.  Thus in $V[G]$, the set $\mathsf{W}$ is the unique
witness to the formula $\chi(w,\omega_2^V)$ which asserts that there exists an increasing sequence
$(\beta_{\xi} \ | \ \xi < \omega_1)$, which is cofinal in $\omega_2^V$, and a sequence $(w_{\xi} \ | \ \xi < \omega_1)$ with each $w_{\xi}$ satisfying $\psi(w_{\xi},\beta_{\xi})$ and $$w = \bigcup_{\xi < \omega_1} w_{\xi}.$$  \noindent
With a $\Sigma_1$ formula involving $\omega_2^V$ as a parameter we can therefore identify $H(\omega_2)^{V}$ as $$H = L_{\omega_2^V}[\mathsf{W}],$$ \noindent and using the definitions of $\mathbb{P}$ and $D$ assert the existence of a filter $G$ on $\mathbb{P}$ which meets
$D(a)$ for every $a \in H$.  Indeed, all of this can be verified by a transitive structure $N$ of a sufficient
fragment of set theory containing $G$ and $H$ which satisfied a formula involving $\omega_2^V$.  Thus in
$V$, by intersecting with the appropriate club, we get a stationary set of $\delta$ such that $$L_{\delta}[\mathsf{W} \upharpoonright \delta] \prec L_{\omega_2}[\mathsf{W}] = H(\omega_2)$$ \noindent and a filter $\bar{G}$ on $\mathbb{P} \cap X$ where $X = L_{\delta}[\mathsf{W} \upharpoonright \delta]$ which has the desired properties.  \end{proof}

\begin{thm} $\mathsf{BMM^{s,++}}$ implies $\mathsf{PD}$ in all generic extensions. \end{thm}
\begin{proof} Recall that $M_{n}^{*}(a)$ is the minimal sound $a$-mouse with active top extender which is closed under the $M_{n}^{\#}$ operation.  We show that $M_{n}^{*}(a)$ exists for every transitive set $a$ by induction on $n < \omega$.  The base case is already accomplished by Theorem 5 so we just assume the induction hypothesis holds for some $n < \omega$.  We need to see that $$M_{n}^{*}(\mathsf{W}) \models NS \mbox{ saturated }.$$
\noindent Suppose toward a contradiction that there is a maximal antichain in $P(\omega_1)/NS$ which belongs to $M_{n}^{*}(\mathsf{W})$ and has size $\omega_2$.  We assume that $\mathcal{A}$ is the least such
in the definable wellordering of $M_{n}^{*}(\mathsf{W})$.  Since $P(\omega_1) \subset N$, the notion of
being a maximal antichain is absolute.
Let $\mathbb{P}$ be the standard poset for sealing $\mathcal{A}$.  Thus if $G \subset \mathbb{P}$ is $V$-generic
then in $V[G]$ there is an enumeration $$\mathcal{A} = \{ A_{\alpha} \ | \ \alpha < \omega_1 \}$$ \noindent
whose diagonal union contains a club $C$.  Inside $V[G]$ let $N$ denote the transitive collapse of
an elementary submodel $X$ of a large enough $H(\theta)$ so that $X$ contains $M_{n}^{*}(\mathsf{W})$ as well as the enumeration of $\mathcal{A}$ and the club $C$.  Then $N$ reflects the relevant properties of these objects mentioned
above.  Under our closure assumptions, there is a formula $\chi$ so
that for any countable and transitive set $a$ and countable structure $M$, $$M = M_{n}^*(a) \Leftrightarrow (H(\omega_1),\in) \models \chi(a,M).$$
\noindent Since $\mathbb{P}$ does not add countable sets we know that

$$H(\omega_1)^{V} = H(\omega_1)^{V[G]} = H(\omega_1)^{L[W]}.$$ \noindent Hence, inside $N$ there is
a continuous sequence of substructures $$\{ X_{\alpha} \ | \ \alpha < \omega_1 \}$$ \noindent of some $H(\kappa)$,
each of which contain $M_{n}^{*}(\mathsf{W})$ so that letting $M_{\xi}$ and $w_{\xi}$ be the
image of $M_{n}^{*}(\mathsf{W})$ and $\mathsf{W}$ respectively under the map which collapses
$X_{\xi}$, the fact that $M_{\xi} = M_{n}^{*}(w_{\xi})$ is certified by the formula $\chi$ and
the structure $H(\omega_1)^{L[W]}$.  Back in $V$ we get a model $\bar{N}$ as above whose
version of $\mathsf{W}$ is $\mathsf{W}_{\delta}$ where $\delta$ is such that $M_{n}^{*}(\mathsf{W}_{\delta})$ is
fully elementary in $M_{n}^{*}(\mathsf{W})$.  It follows that the version of $M_{n}^{*}(\mathsf{W}_{\delta})$ that
$\bar{N}$ sees is the true version, since it collapsed correctly on a club.  Moreover, $\bar{N}$ sees that the least antichain $\mathcal{A}_{\delta}$ of $M_{n}^{*}(\mathsf{W}_{\delta})$ is sealed.  Since this antichain is a subset of $\mathcal{A}$ this gives
a contradiction.  Now, under these conditions, the argument of Lemmas 16 and 17 of \cite{SZ} immediately give closure of $P(\omega_1)$ under the $M_{n+1}^{\#}$ operation, and we turn toward
the lifting portion of the induction step.  This is modeled on 10.108 of \cite{W} which shows that
$\mathsf{BMM}^{++}$ lifts closure under sharps from $P(\omega_1)$ to all of $V$.  We claim that
$M_{n+1}^{\#}(a)$ exists for every set $a$.  Otherwise we may pass to $V[g]$ where $g \subset Col(\omega,\kappa)$ is $V$-generic for a sufficiently large $\kappa$ and find a subset $b$ of $\omega_1$, a term $t$, and stationary sets $S,T$ such that $$\beta \in S \Rightarrow t \in M_{n+1}^{\#}(b \cap \beta)$$ \noindent and
$$\beta \in T \Rightarrow t \notin M_{n+1}^{\#}(b \cap \beta).$$  \noindent Using the same trick
as above to certify each $M_{n+1}^{\#}(b \cap \beta)$, we find in sets $\bar{b},\bar{S},\bar{T}$ back in
$V$ with the property above.  This contradicts the existence of $M_{n+1}^{\#}(\bar{b})$.  Repeating the
two arguments finishes the proof.  \end{proof}

\begin{thm} $\mathsf{BMM^s}(c)$ fails in the $\mathbb{P}_{max}$ model for $\mathsf{BMM}$.
\end{thm}

\begin{proof} Assume otherwise.  Thus we have $\mathsf{BMM^s}(c)$, $\mathsf{BMM}$
and NS saturated at our disposal.  We will prove $\mathsf{PD}$ from these assumptions and the proof will yield the
desired contradiction.   We first claim that if $S \subset [\omega_2]^{\omega}$ is
projective stationary and first order definable over the structure $(H(\omega_2),\in)$.  Then $S \cap [\delta]^{\omega}$ \noindent contains a club in $[\delta]^{\omega}$ for a stationary set of $\delta < \omega_2$.  This principle, denote $\mathsf{DSRP^s}(\omega_2)$, can be deduced from Lemma 9 as follows.  Suppose $S$ is such a set and let $\mathbb{P}$ be the standard poset for shooting a club through $S$.  Thus elements of $\mathbb{P}$ are countable continuous increasing sequences $$p = (\sigma_{\xi} \ | \ \xi \leq \gamma)$$ \noindent from $S$.  For $\alpha < \omega_2$ let $D(\alpha)$ denote the set of conditions $p$ as above for which $\alpha \in \sigma_{\xi}$ for some $\xi$ in the domain of $p$.  Lemma 9
gives the desired stationary set of club reflection points.  We now claim that if $S,T \subset [\omega_2]^{\omega}$ are stationary and first order definable over the structure $(H(\omega_2),\in)$.  Then $$S \cap [\delta]^{\omega} \mbox{ and } T \cap [\delta]^{\omega}$$ \noindent are both stationary in $[\delta]^{\omega}$ for a stationary set of $\delta < \omega_2$.  Given such
a pair $S,T$ it follows from NS saturated that there are stationary sets $A_{S},A_T \subset \omega_1$ such that
$S$ is $A_{S}$-projective stationary, $T$ is $A_{T}$-projective stationary, and $$A_{S} \cap A_{T} = \emptyset.$$
\noindent The set
$$P(S,T) = \{ \sigma \ | \ (\sigma \cap \omega_1 \in A_{S} \rightarrow \sigma \in S) \wedge (\sigma \cap \omega_1 \in A_{T} \rightarrow \sigma \in T) \}$$ \noindent is projective stationary and so reflects to a
club in $[\delta]^{\omega}$ for a stationary set of $\delta < \omega_2$ by $\mathsf{DSRP^s}(\omega_2)$, and this 
proves the claim.  Woodin's proof of $\mathsf{PD}$ from $\mathsf{MM}(c)$ only uses
NS saturated and the simultaneous reflection principle $\mathsf{WRP}_{(2)}(\omega_2)$.  The definable version of
this principle that we now have at our disposal suffices with the caveat that one can only show that $M_{n}^{*}(\mathsf{W})$ exists by induction on $n < \omega$, as opposed to closure of $P(\omega_2)$ under the $M_{n}^{*}$ operation.  This however, is enough to implement the argument, and we refer the reader to \cite{SZ} for more details.  Since $M_{1}^{*}(\mathsf{W})$ does not exist in the $\mathsf{BMM}$ model we get the desired contradiction.  \end{proof}

\begin{thm} Let $N$ be the minimal inner model containing $\mathbb{R}$ and closed under
the $M_{1}^{\#}$ operation, and assume $N \models \mathsf{AD}$.  Then
$$N[G] \models \mathsf{BMM^{s_0,++}}$$
\noindent whenever $G \subset \mathbb{P}_{max}$ is $N$-generic.
\end{thm}

\begin{proof}  Suppose $G \subset \mathbb{P}_{max}$ is
$N$-generic and $\mathbb{Q}$ is a poset in $N[G]$ such that $$\Vdash_{\mathbb{Q}}^{N[G]} \phi(\omega_2^V,a^{*})$$
\noindent holds where $\phi(x,a)$ is $\Sigma_1$ in the appropriate language with parameter
$a^* \subset \omega_1$.  We may assume that $a^* = a_{G}$.  We also assume that
$$\Vdash_{\mathbb{Q}}^{N[G]} \mbox{cf}(\omega_2^{V}) = \omega.$$

\noindent Fix a condition $\mathcal{M}_{0} = ((M_0,I_0),a_0) \in G$ which forces this as well as that $\dot{C}$ is a club subset of $\omega_2$.  Note that $$H(\omega_1)^{N} = H(\omega_2)^{N[G]}.$$  \noindent Working in $N[G]$ we are going to produce a condition $\mathcal{M}_1$ below $\mathcal{M}_{0}$ so that $\mathcal{M}_1 \in N$ and
 $$\mathcal{M}_1 \Vdash^{N}_{\mathbb{P}_{\max}} \exists \gamma \ \gamma \in \dot{C} \wedge \phi(\gamma,\dot{a}_{G}).$$
\noindent This will prove the theorem.   First let us introduce some notation.  We think of $x^{\#}$ for a real $x$ as an $x$-mouse $(L_{\alpha}[x],\mu)$.  Let $\kappa$ be the critical point of the measure $\mu$, and $j$ the map obtained by iterating the measure $\omega_1$ times.  We say that a pair $c = (x^{\#},\beta)$ with $\kappa < \beta < \alpha$ is a {\em code} for an ordinal $\gamma$ if $j(\beta) = \gamma$.  We let $\gamma_{c}$ denote the ordinal just described.  In our present situation, every ordinal less
 than $\omega_2^{V}$ has a code because $$u_2 = \undertilde{\delta}^1_2 = \omega_2$$ \noindent in $N[G]$.  Now, let $H \subset \mathbb{Q}$ be $N[G]$
generic.  Using our closure hypothesis, we can create a condition $\mathcal{M} = ((M,I),a^*)$ in a sufficiently large collapse over $N[G][H]$ with an ordinal $\delta \in M$ satisfying the following conditions.  Note that $\mathcal{M}_{0}$ is iterable in all generic extensions as $N[G]$ is sufficiently correct.

\begin{enumerate}

\item $((M,I),a^*) < ((M_0,I_0),a_0)$

\item $M \models \phi(\delta,a^*)$

\item $M \models f:\omega \rightarrow \delta$ is a cofinal

\item for every $n < \omega$ there is a condition $\mathcal{P}_{n} = ((P_n,J_{n}),b_{n})$ which
is greater than $((M,I),a^*)$ and a code $c(n)$ such that

\begin{enumerate}

\item $\mathcal{P}_{n}$ and $c(n)$ belong to $M$ and $\mathcal{M} < \mathcal{P}_{n}$

\item $\mathcal{P}_{n} \Vdash^{N}_{\mathbb{P}_{\max}} \gamma_{c(n)} \in \dot{C}$

\item $M \models f(n) < \gamma_{c(n)} < \delta$

\end{enumerate}

\end{enumerate}

\noindent  By $<$ we mean of course the $\mathbb{P}_{max}$ ordering.  To construct the condition let $f: \omega \rightarrow \delta$ be any cofinal map where $\delta = \omega_2^{V}$ and let $\theta$ be sufficiently large.  Note that $$H(\theta) \models \phi(\delta,a^*).$$  \noindent Let $E$ be set of ordinals so that $H(\theta) \in L[E]$.  Let $$Y = M_{1}^{\#}(E)$$ \noindent
and let $g$ be generic over $N[G][H]$ so that $Y$ is countable in $N[G][H][g]$.  Let $\hat{g} \in N[G][H][g]$ be $Y$ generic for making NS presaturated and then forcing $\mathsf{MA}$, and let $$\mathcal{M} = ((Y[\hat{g}],NS^{Y[\hat{g}]}),a^*).$$
\noindent We claim that $\mathcal{M}$ satisfies the conditions above.  Since $P(\omega_1)^{N[G][H]} \subset Y$ and $H$ is generic over $N[G]$ for stationary set preserving forcing we know that $$NS^{Y} \cap N[G] = NS^{N[G][H]} \cap N[G] = NS^{N[G]},$$ \noindent and since $\hat{g}$ preserves stationary sets we have $$NS^{Y[\hat{g}]} \cap N[G] = NS^{N[G]}.$$  \noindent It follows that (1) holds as witnessed by the iteration of $\mathcal{M}_{0}$ determined by the generic
$G$.  The next two conditions hold as they are upward absolute.  For $n < \omega$ there will be a condition $\mathcal{P}_{n} \in G \cap H(\theta)^{N[G][H]}$ and a code $c(n)$ with the properties above because
$\dot{C}_{G}$ is cofinal in $\delta = \omega_2^{V}$, and by the reasoning used to establish (1).  Now let us go back to $N[G]$.  Let $X$ be an elementary submodel of a large enough rank initial segment of $N[G]$
which contains everything relevant and let $\pi:X \rightarrow \mathcal{N}$ denote the transitivization map.
Let $\bar{H} \subset \pi(\mathbb{Q})$ be $\mathcal{N}$ generic for the collapse of $\mathbb{Q}$ and let $\bar{g}$
be $\mathcal{N}[\bar{H}]$ generic for the sufficiently large collapse in the sense of $\mathcal{N}$.  Then $\mathcal{N}[\bar{H}][\bar{g}]$ thinks there is a condition $$\mathcal{M}_{1} = ((M_1,I_1),a_1)$$
\noindent which satisfies the conditions above.  Since $\mathcal{N}$ is closed under sharps, this condition is truly iterable in $N[G]$.  Of course $a_1 = a^* \cap \mathcal{N}$ but the rest of the properties are upward absolute.  The second clause of condition (4) holds by elementarity of $\pi$ and the fact that the conditions $\mathcal{P}_{n}$ are countable.  Moreover, this condition is in the ground model as $\mathbb{P}_{max}$ does
not add reals, and has the properties there as well.  Let us check that $$\mathcal{M}_1 \Vdash^{N}_{\mathbb{P}_{\max}} \exists \gamma \ \gamma \in \dot{C} \wedge \phi(\gamma,\dot{a}_{G}).$$  \noindent Let $G \subset \mathbb{P}_{max}$
be $N$ generic below $\mathcal{M}_1$ and let $$j:M_1 \rightarrow M^*$$ \noindent be the iteration determined
by $G$, and let $\delta$ and $f$ be as in the conditions enumerated above.  Thus
$\mathcal{M}^{*} \models \phi(j(\delta),a_{G})$ and hence $$(H(\omega_2),\in,NS) \models \phi(j(\delta),a_G).$$
Let $C = \dot{C}_{G}$.  For each $n < \omega$ we have $\mathcal{M}_1 < \mathcal{P}_{n}$ and so $\mathcal{P}_{n} \in G$ as well.  Thus $\gamma_{c(n)} \in C$ for each $n < \omega$.  Now, we may assume that the sequence
$$(\mathcal{P}_{n} \ | \ n < \omega)$$ \noindent is an element of $M_1$ although this is not necessary.  Thus

$$M_1 \models \delta = \bigcup_{n<\omega} \gamma_{c(n)}$$ \noindent where each $\gamma_{c(n)}$ is computed
in $M_1$, and so
$$M^* \models j(\delta) = \bigcup_{n<\omega} \gamma_{c(n)}$$ \noindent and we conclude
that $j(\delta) \in C$ as desired.\footnote{If $\omega_2^V$ were not countably cofinal in $N[G][H]$
we could choose $f: \omega_1 \rightarrow \delta$ to be a bijection and use conditions $\mathcal{P}_{\xi}$ for $\xi < \omega_1$ as above.  The problem occurs at the end of the argument as $\delta$ is not a continuity point of the embedding $j$.} \end{proof}

The proof given above, with the extra
moves required for the $s_0$ clause suitable excised, constitutes a reorganization of the proof of the
following equivalent formulation of the consistency result for $\mathsf{BMM}$ from \cite{W}.

\begin{quote}  Assume $(*)$ and that $M_{1}^{\#}(X)$ exists for every set.  Suppose
$N$ is an inner model of $\mathsf{ZFC}$ containing $P(\omega_1)$ and closed under
the $M_{1}^{\#}$ operation.  Then $N \models \mathsf{BMM}^{++}$. (10.99 of \cite{W})\end{quote}

However, even though $\mathsf{BMM^{s_0}}$ holds in the $\mathsf{BMM}$ model, it is
not a consequence of $(*)$ together with global closure under the $M_{1}^{\#}$ operation.  This phenomenon
is well precedented in \cite{W}, for example in the case of the saturation of the nonstationary ideal, which
holds in the $\mathbb{P}_{max}$ extension of $L(\mathbb{R})$ but is not a consequence of the $\mathbb{P}_{max}$
axiom $(*)$.  Recall that $\tilde{T}$, for a set $T \subset \omega_1$, is the set of $\alpha < \omega_2$ for which there is a club of $\sigma \in [\alpha]^{\omega}$ with the order type of $\sigma$ in $T$.  Theorem 5.8 of \cite{W},
which was used in the proof of Theorem 7, shows that under $\mathsf{MM}$ the set
$$\tilde{T}^{0} = \{ \alpha \in \tilde{T} \ | \  cf(\alpha) = \omega \},$$ \noindent is stationary for
every stationary set $T \subset \omega_1$.  It is straightforward to check that $\mathsf{BMM^{s_0}}$ together with the saturation of the nonstationary ideal suffice for this result.  Arguments of Larson from \cite{L} can be used to show that the poset $\mathbb{P}$
for shooting an $\omega$-club through $\tilde{T}^{0}$ over the $\mathsf{BMM}$ model does not add new subsets of $\omega_1$ and hence preserves $(*)$
together with global closure under the $M_{1}^{\#}$ operation.  Assuming $T$ is costationary this
yields the desired separation.  Moreover, his arguments show that $\mathsf{MM}^{++}(c)$ could
be preserved as well if the ground model were taken to be the richer $\mathbb{P}_{max}$ model for
$\mathsf{MM}(c)$ together with $\mathsf{BMM}$.  For a proof that $\mathsf{BMM^{s_0}}$ is not
implied by $\mathsf{BMM}^{++}$ the interested reader could just check that $\mathbb{P}$ does not add $\omega_1$ sequences under the assumption that $\mathsf{MM}$ holds in the ground model.

Finally, we prove another separation result which does not seem to involve the consequences of
$\mathsf{BMM^s}$ for projective stationary sets.  It involves rather the concept of
a {\em disjoint club sequence on} $\omega_2$, which is a sequence $$(C_{\alpha} \ | \ \alpha \in A)$$ \noindent of pairwise disjoint sets, with each $C_{\alpha}$ a club subset of $[\alpha]^{\omega}$ and $A$ a stationary subset of $\omega_2$ consisting a ordinals of uncountable cofinality.  This is an invention of Krieger from \cite{K} who derives one from $\mathsf{MM}(c)$.

\begin{thm}  $\mathsf{BMM^s}$ implies the existence of a disjoint club sequence on $\omega_2$.
\end{thm}

\begin{proof}  Let us fix a canonical way of coding sets like
$C_{\alpha}$ above as subsets of $\omega_1$.  Define
$A_{\mathsf{W}} \subset \omega_2$ and $\vec{C} = \{ C_{\alpha} \ | \ \alpha \in A_{\mathsf{W}}\}$ by induction as follows.  Given $\vec{C} \upharpoonright \alpha$ and $A_{\mathsf{W}} \cap \alpha$, for an ordinal
$\alpha$ of uncountable cofinality, if $$C = \bigcup_{\beta < \alpha} C_{\beta} \subset [\alpha]^{\omega}$$ \noindent is nonstationary in $[\alpha]^{\omega}$ then let $C_{\alpha}$ be a club disjoint from $C$ with the earliest index according
to $\mathsf{W}$, and put $\alpha$ in $A_{\mathsf{W}}$.  We need to see that $A_{\mathsf{W}}$ is stationary.  Theorem 4.4 of \cite{K} shows that there is a stationary set preserving
notion of forcing $\mathbb{P}$ such that whenever $G \subset \mathbb{P}$ is $V$-generic $\omega_2^V$ has uncountable cofinality and there is in $V[G]$ a club $C$ in $[\omega_2^V]^{\omega}$ which is disjoint from
$$D = \bigcup_{\beta < \omega_2^V} C_{\beta} \subset [\alpha]^{\omega}.$$  \noindent The key point
is that for any $\sigma \in C$ it can be verifies that $\sigma \notin D$ by consulting
$H = L_{\omega_2^{V}+1}[\mathsf{W}]$ of which $\vec{C}$ is an element.  The existence of club $C$ and the structure
$H$ witnessing that $C \cap D = \emptyset$ is a $\Sigma_1$ property of $\omega_2^V$
so we get a stationary set of witnesses $\delta < \omega_2$ in $V$, each of which such that
$$L_{\delta +1}[\mathsf{W}] \prec L_{\omega_2 +1}[\mathsf{W}],$$
\noindent and each of these ordinals must therefore belong to $A_{\mathsf{W}}$ as desired.  \end{proof}

\noindent The argument of 3.4 of \cite{K} which shows that $$A \cup \{ \gamma < \omega_2 \ | \ cf(\gamma) = \omega \}$$ \noindent does not contain a club whenever $A$ indexes a disjoint club sequence is used to show
that any disjoint club sequence can be killed with a forcing that leaves $H(\omega_2)$ undisturbed.  We are sure this would be known to the authors of \cite{K} but we prove it here so we can observe that in the extension the
analogue of $\mathsf{BSRP^s}(\omega_2)$ from Theorem 11, which we denote by $\mathsf{DSRP^s}(\omega_2)$, persists while the set $A_{\mathsf{W}}$ becomes nonstationary so that $\mathsf{BMM^{s}}(c)$ fails.

\begin{thm}  Assume $\mathsf{MM}$.  Then there is a forcing notion $\mathbb{P}$ of size $\omega_2$ such that whenever $G \subset \mathbb{P}$ is $V$-generic, $$V[G] \models \mathsf{BMM}^{++} \wedge \mathsf{DSRP^s}(\omega_2)
\wedge \neg \mathsf{BMM^{s}}(c).$$
\end{thm}

\begin{proof}  Let $\{ C_{\alpha} \ | \ \alpha \in A_{\mathsf{W}}\}$ be the set produced in Theorem 7.  Let $\mathbb{P}$ consist of closed
subsets of $\omega_2 \setminus \A$ of size $\omega_1$, ordered by end extension.
We claim that forcing with $\mathbb{P}$ does not introduce new subsets of $\omega_1$.  Note that $\mathbb{P}$ is
$\sigma$-closed.  Let
$\tau$ be a $\mathbb{P}$ term which is forced by a condition $p$ to be a subset of $\omega_1$.  Fix
a large enough $\theta$, and consider sequences $$((N_{\gamma},p_{\gamma},s_{\gamma}) \ | \ \gamma < \omega_1)$$ \noindent satisfying the following conditions:

\begin{enumerate}

\item $N_{\gamma} \prec H(\theta)$ and $N_{\gamma} \in N_{\gamma+1}$

\item $(N_{\gamma} \ | \ \gamma < \omega_1)$ is increasing and continuous

\item $p_{0} = p$ and each $p_{\gamma} \in \mathbb{P} \cap N_{\gamma}$

\item $( p_{\gamma} \ | \ \gamma < \omega_1 )$ is $<_{\mathbb{P}}$-decreasing

\item $p_{\gamma} \Vdash_{\mathbb{P}} \tau \cap \gamma = s_{\gamma}$

\item $\{ N_{\gamma} \cap \omega_2 \ | \ \gamma < \omega_1 \}$ is club in $[\alpha]^{\omega}$ for some $\alpha< \omega_2$.

\end{enumerate}

\noindent We need find such a sequence with $\alpha \notin{\A}$, for then
$$q = (\bigcup_{\gamma < \omega_1} p_{\gamma}) \cup \{\alpha \} \in \mathbb{P}$$ \noindent
and $$q \Vdash_{\mathbb{P}} \tau = f \mbox{ where } f = \bigcup_{\gamma < \omega_1} f_{\gamma}.$$
\noindent Define $B$ to be the set of $\alpha$ for which there exists a sequence as above.  It is
easy to see that $B$ is stationary.  So suppose toward a contradiction that $B \subset \A$ and
for $\alpha \in B$ and let $(N^{\alpha}_{\gamma} \ | \ \gamma < \omega_1)$ be the sequence as above.
As in 3.4 of \cite{K}, let $c_{\alpha} \subset \omega_1$ be club so that
$$\{ N^{\alpha}_{\gamma} \cap \omega_2 \ | \ \gamma \in c_{\alpha} \}$$ \noindent is club in
$[\alpha]^{\omega}$ and contained in $C_{\alpha}$.  Let $i_{\alpha}$ be the minimum element of $c_{\alpha}$
and let $d_{\alpha}= c_{\alpha} \setminus \{i_{\alpha}\}$.  Let $$S = \{ N^{\alpha}_{\gamma} \cap H(\omega_2) \ | \ \alpha \in B \wedge \gamma \in d_{\alpha} \}.$$  Then $S$ is stationary in $[H(\omega_2)]^{\omega}$ and
by pressing down we get $\alpha < \beta$ such that $$N^{\alpha}_{i_{\alpha}} \cap \omega_2 = N^{\beta}_{i_{\beta}} \cap \omega_2,$$ a contradiction as $C_{\alpha} \cap C_{\beta}$ is empty.  Now let $G \subset \mathbb{P}$ be
$V$-generic.  We have that $$H(\omega_2)^{V[G]} = H(\omega_2)^{V}$$ \noindent and
so $V[G] \models \mathsf{BMM^{++}}$ and $$(\A)^{V[G]} = (\A)^{V[G]}$$ \noindent so $\mathsf{BMM^{s,++}}(c)$
must fail in $V[G]$.  Now fix a projective stationary set $S \subset [\omega_2]^{\omega}$ which
is first order definable over $H(\omega_2)^{V[G]}$.  Thus $S$ is projective stationary in $V$.  Let $\dot{C}$
be a term for a club subset of $\omega_2$.  Note that the proof above shows that $B^* = B \setminus \A$ is
stationary.  We may assume that we have required that the condition produced at the end forces that
$\alpha \in \dot{C}$.  Using $\mathsf{MM}$ we can then find such an $\alpha \in B^*$ which is
a club reflection point for $S$.  Thus $\mathsf{DSRP^s}(\omega_2)$ continues to hold.  \end{proof}

\end{document}